\documentclass[12pt,twoside]{article}
\usepackage{amsmath,amssymb}
\newtheorem{Def}{Definition}[section]

\newtheorem{lem}[Def]{Lemma}
\newtheorem{stw}[Def]{Proposition}
\newtheorem{wn}[Def]{Corollary}

\newcommand{\NN}{{\rm I\!N}}
\newcommand{\RR}{{\rm I\!R}}
\newcommand{\CC}{{\mathbb C}}

\newcommand{\be}{\begin{equation}}
\newcommand{\ee}{\end{equation}}
\newcommand{\bsm}{\left(\begin{smallmatrix}}
\newcommand{\esm}{\end{smallmatrix}\right)}
\newcommand{\dow}{{\bf Proof: }}

\newcommand{\qed}{\hfill{Q.E.D.}\vspace{3.6mm}}
\newcommand{\where}{,\hspace{7mm}}

\newcommand{\Int}{\int\limits}
\newcommand{\Sp}{{\rm Sp\,}}
\newcommand{\G}{{\rm Graph}}
\renewcommand{\Im}{{\rm Im\,}}

\newcommand{\F}{{\cal F}}
\newcommand{\la}{\langle}
\newcommand{\ra}{\rangle}

\begin{document}

\title{A remark on the spectrum of the analytic generator}
\author{Piotr Miko\l{}aj So\l{}tan\thanks{Partially
supported by Komitet Bada\'{n} Naukowych, grant No 2 P0A3 030 14}\\
Department of Mathematical Methods in Physics\\
Faculty of Physics, University of Warsaw\\
Ho\.{z}a 74, 00-682 Warsaw\\
Poland\\
pmsoltan@polbox.com}
\date{}
\maketitle
\begin{abstract}
\noindent
It is shown that although the spectrum of the analytic generator of a
one--parameter group of isometries of a Banach space may be equal to $\CC$
(cf \cite{vd} and \cite{elzs}), a simple operation of ampliating the analytic
generator onto its graph locates its spectrum in $\RR_+$.
\end{abstract}

\section{Introduction}

Let $(X,\F)$ be a dual pair of Banach spaces with the pairing denoted by
$\la\cdot,\cdot\ra$. The symbol $\sigma(X,\F)$ will denote the weak topology
on $X$ given by the pairing with $\F$ as well as the product of such
topologies on $X\times X$. We shall say that {\it the pair $(X,\F)$ has the
Krein property} if the \hbox{$\sigma(X,\F)$-closed} convex hull of any
\hbox{$\sigma(X,\F)$-compact} set in $X$ is again
\hbox{$\sigma(X,\F)$-compact}. Throughout the paper we shall assume that both
$(X,\F)$ and $(\F,X)$ have the Krein property. Let $B_\F(X)\subset B(X)$ be
the subspace of \hbox{$\sigma(X,\F)$-continuous} linear maps of $X$ into
itself. Let $U=\{U_t\}_{t\in\RR}$ be a \hbox{$\sigma(X,\F)$-continuous},
one--parameter group of isometries in $B_\F(X)$. We shall denote by $X_\infty$
the \hbox{$\sigma(X,\F)$-dense} subspace of $X$ consisting of entire analytic
elements for $U$. For any $z\in\CC$ the operator $U_z$ is defined in the
following way:
\[
\left(\begin{pmatrix}x\\y\end{pmatrix}\in\G(U_z)\right)
\Longleftrightarrow
\left(
\begin{array}{c}
\mbox{There exists a $\sigma(X,\F)$-continuous}\\
\mbox{function $F_x$ defined on the strip}\\
\{w\colon \Im{w}\,\Im{z}\geq 0,\ |\Im{w}|\leq|\Im{z}|\},\\
\mbox{with values in $X$, holomophic inside}\\
\mbox{the strip, such that for all }t\in\RR\\
\mbox{$F_x(t)=U_tx$ and }F_x(z)=y\\
\end{array}
\right).
\]
It is known (cf \cite{cizs}, \cite{elzs}) that $U_z$ is a well-defined,
\hbox{$\sigma(X,\F)$-closed} linear operator in $X$ and that
$U_{z_1}U_{z_2}=U_{z_1+z_2}$ for all $z_1,z_2\in\CC$. The {\it analytic
generator of the one--parameter group $U$} is the operator $U_i$.

\noindent
At first glance it seems that the spectrum of the analityc generator of $U$
should be contained in $\RR_+$. In fact the analytic generator of a strongly
continuous group of unitaries on a Hilbert space is a positive, self-adjoint
operator. However as it was first shown in \cite{vd} it may happen that
$\Sp U_i=\CC$. Furthermore it turns out that we have either
$\Sp U_i\subset\RR_+$ or $\Sp U_i=\CC$ (cf \cite{elzs}).

\noindent
We shall use some results on integration of vector valued functions (cf
\cite{arv} sect. 1). Given a locally compact space $\Omega$, a complex,
regular, Borel measure of finite variation $\nu$ on $\Omega$ and a
\hbox{$\sigma(X,\F)$-continuous}, norm-bounded function
$\Omega\ni\omega\mapsto x(\omega)\in X$, the Krein property of $(X,\F)$ implies
that there is a unique $y\in X$ with
\[
\la y,\phi\ra=\Int_\Omega \la x(\omega),\phi\ra d\nu(\omega)\where \phi\in\F.
\]
As usual we shall write $y=\Int_\Omega x(\omega)d\nu(\omega)$. Similarly for
any complex, regular, Borel measure of finite variation $\nu$ on $\RR$ the
operator
\[
X\ni x\longmapsto \Int_\RR U_txd\nu(t)\in X
\]
will be denoted by $\Int_\RR U_td\nu(t)$. This operator is bounded and, in
fact, thanks to the Krein property of $(\F,X)$ it is
\hbox{$\sigma(X,\F)$-continuous}.

\noindent
The pairing of $X$ and $\F$ gives rise to the operation of transposition
defined on the set of \hbox{$\sigma(X,\F)$-densly} defined,
\hbox{$\sigma(X,\F)$-closed} operators on $X$. This operation shall be denoted
by $S\mapsto S^{\rm T}$.

\begin{stw}
\label{komut}
Let $(X,\F)$ be a dual pair of Banach spaces such that both $(X,\F)$ and
$(\F,X)$ have the Krein property. Let $U=\{U_t\}_{t\in\RR}$ be a
\hbox{$\sigma(X,\F)$-continuous} group of isometries in $B_\F(X)$ and let
$S$ be a \hbox{$\sigma(X,\F)$-densly} defined, \hbox{$\sigma(X,\F)$-closed}
operator such that $U_tS=SU_t$ for all $t\in\RR$. Let $\nu$ be a complex,
regular Borel measure of finite variation on $\RR$ and denote
$A=\Int_\RR U_td\nu(t)$. Then $AS\subset SA$.
\end{stw}
\dow
For any $\phi\in D(S^{\rm T})=\bigl(\mbox{the domain of }S^{\rm T}\bigr)$ and
$x\in D(S)$ we
have
\begin{eqnarray*}
\la ASx,\phi\ra\hspace{.22cm}
=&\Int_{\RR}\la U_tSx,\phi\ra d\nu(t)=\Int_{\RR}\la SU_tx,\phi\ra d\nu(t)&\\
=&\hspace{-.84cm}\Int_{\RR}\la U_tx,S^{\rm T}\phi\ra d\nu(t)=\la Ax,S^{\rm
T}\phi\ra&
\end{eqnarray*}
which means simply that $Ax\in D(S^{\rm TT})=D(S)$ and $SAx=S^{\rm TT}Ax=ASx$.
Thus $AS\subset SA$.
\qed

\noindent
In the next section we shall see that a simple operation on $U_i$ can squeeze
its spectrum into $\RR_+$. Let $\Delta$ be the ampliation of $U_i$ onto
its graph: for $\bsm x\\x'\esm,\bsm y\\y'\esm\in\G(U_i)$ we define
\[
\left(
\begin{array}{c}
\bsm x\\x'\esm \in D(\Delta),\\
\bsm y\\y'\esm = \Delta\bsm x\\x'\esm
\end{array}
\right)
\Longleftrightarrow
\left(
\begin{array}{c}
x,x'\in D(U_i),\\
y=U_ix,\ y'=U_ix'
\end{array}
\right).
\]

\section{The spectrum of $\Delta$}

Let $\mu\in\CC\setminus\RR_-$ be a parameter. Define an integrable function
$F_\mu=F_\mu(t)$ on $\RR$ by
\[
F_\mu(t)=\frac{1}{2\pi}\Int_{-\infty}^{+\infty}\frac{e^{E(1+it)}}{(e^E+\mu)^2}dE
=\frac{t(-\mu)^{it-1}}{e^{-2\pi t}-1}=\frac{t\mu^{it-1}}{e^{\pi t}-e^{-\pi t}}.
\]\sloppy
For any fixed $\mu$ the function $F_\mu$ has holomorphic continuation onto
the set \hbox{$\CC\setminus\{\pm ni\colon n=1,2,3,\ldots\}$} and for
$0\neq z$ in that region we have
\begin{eqnarray}
&F_\mu(t-2i)+2\mu F_\mu(t-i)+\mu^2F_\mu(t)=0,\label{wl1}&\\
&\mu F_\mu(z)+F_\mu(z-i)=i\frac{\mu^{iz}}{e^{\pi z}-e^{-\pi z}}.\label{wl2}&
\end{eqnarray}
Define now a linear operator $Q_\mu\colon X\rightarrow X$
\begin{equation}
\label{Q}
Q_\mu=\Int_{-\infty}^{+\infty} F_\mu(t)U_tdt\in B_\F(X).
\end{equation}
Take $x\in X_\infty$. We have
\begin{equation}
\label{tu}
\left.
\parbox{11.9cm}
{
\vspace{-.5cm}
\begin{eqnarray*}
&\hspace{-6.26cm}(Q_\mu U_{2i}+2\mu Q_\mu U_i+\mu^2Q_\mu)x&\\
&\hspace{-.32cm}=\Int_{-\infty}^{+\infty}F_\mu(t)U_{t+2i}xdt+2\mu\Int_{-\infty}^
{+\infty}F_\mu(t)U_{t+i}xdt
+\mu^2\Int_{-\infty}^{+\infty}F_\mu(t)U_txdt&\\
&\hspace{-.82cm}=\Int_{C_0}F_\mu(t)U_{t+2i}xdt+2\mu\Int_{C_0}F_\mu(t)U_{t+i}xdt+
\mu^2
                \Int_{C_0}F_\mu(t)U_txdt&\\
&=\Int_{C_2}F_\mu(t-2i)U_txdt+2\mu\Int_{C_1}F_\mu(t-i)U_txdt+
\mu^2\Int_{C_0}F_\mu(t)U_txdt,&
\end{eqnarray*}
\vspace{-.42cm}
}
\right\}
\end{equation}
where $C_0,C_1$ and $C_2$ are oriented curves in $\CC$, as shown on the
following figure:

\noindent
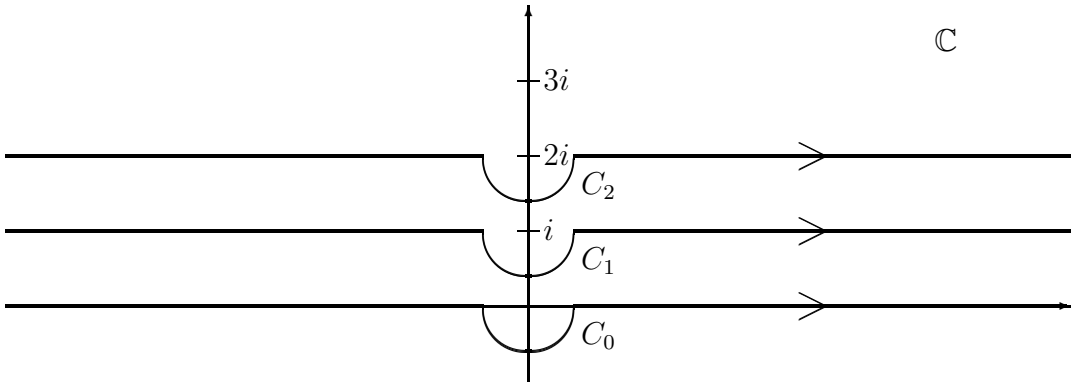
\setlength{\unitlength}{.1mm}
\begin{picture}(1520,600)(0,-50)
\put(45,100){\vector(1,0){1415}}
\put(45,101){\line(1,0){635}}
\put(45,99){\line(1,0){635}}
\put(800,101){\line(1,0){655}}
\put(800,99){\line(1,0){655}}
\put(45,200){\line(1,0){635}}
\put(45,201){\line(1,0){635}}
\put(45,199){\line(1,0){635}}
\put(800,200){\line(1,0){660}}
\put(800,201){\line(1,0){660}}
\put(800,199){\line(1,0){660}}
\put(45,300){\line(1,0){635}}
\put(45,301){\line(1,0){635}}
\put(45,299){\line(1,0){635}}
\put(800,300){\line(1,0){660}}
\put(800,301){\line(1,0){660}}
\put(800,299){\line(1,0){660}}
\put(740,0){\vector(0,1){500}}
\put(725,400){\line(1,0){30}}
\put(760,388){3$i$}
\put(725,300){\line(1,0){30}}
\put(760,288){2$i$}
\put(725,200){\line(1,0){30}}
\put(760,188){$i$}
\put(740,100){\oval(120,120)[b]}
\put(740,200){\oval(120,120)[b]}
\put(740,300){\oval(120,120)[b]}
\put(740,101){\oval(120,120)[b]}
\put(740,201){\oval(120,120)[b]}
\put(740,301){\oval(120,120)[b]}
\put(740,99){\oval(120,120)[b]}
\put(740,199){\oval(120,120)[b]}
\put(740,299){\oval(120,120)[b]}
\put(810,50){$C_0$}
\put(810,150){$C_1$}
\put(810,250){$C_2$}
\put(1280,440){$\CC$}
\put(1100,82){\line(2,1){36}}
\put(1100,118){\line(2,-1){36}}
\put(1100,182){\line(2,1){36}}
\put(1100,218){\line(2,-1){36}}
\put(1100,282){\line(2,1){36}}
\put(1100,318){\line(2,-1){36}}
\end{picture}

\noindent
Now using (\ref{wl1}) we can subtract 0 from both sides of (\ref{tu}) and
taking into account the holomorphy properties of the integrated functions
we obtain
\begin{eqnarray*}
&\hspace{-6.14cm}(Q_\mu U_{2i}+2\mu Q_\mu U_i+\mu^2Q_\mu)x&\\
&=\Int_{C_2}F_\mu(t-2i)U_txdt+2\mu\Int_{C_1}F_\mu(t-i)U_txdt
+\mu^2\Int_{C_0}F_\mu(t)U_txdt&\\
&\hspace{1.38cm}-\left(\Int_{C_2}F_\mu(t-2i)U_txdt+2\mu\Int_{C_2}F_\mu(t-i)U_txd
t+\mu^2\Int_{C_2}F_\mu(t)
U_txdt\right)&\\
&\hspace{-2.56cm}=2\mu\Int_{C_1-C_2}F_\mu(t-i)U_txdt+\mu^2\Int_{C_0-C_2}F_\mu(t)
U_txdt&\\
&\hspace{-3.82cm}=2\mu\Int_{\Gamma}F_\mu(t-i)U_txdt+\mu^2\Int_{\Gamma}F_\mu(t)U_
txdt,&
\end{eqnarray*}
where $\Gamma$ is the oriented curve in $\CC$, as shown below:

\noindent
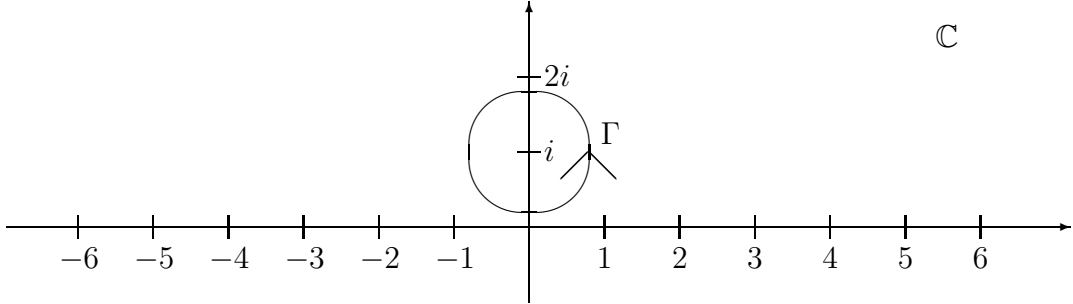
\setlength{\unitlength}{.1mm}
\begin{picture}(1520,500)(0,-50)
\put(45,100){\vector(1,0){1415}}
\put(740,0){\vector(0,1){400}}
\put(725,300){\line(1,0){30}}
\put(760,288){2$i$}
\put(725,200){\line(1,0){30}}
\put(760,188){$i$}
\put(1280,340){$\CC$}
\put(740,200){\oval(160,160)}
\put(835,210){$\Gamma$}
\put(782,164){\line(1,1){36}}
\put(856,164){\line(-1,1){36}}
\put(140,85){\line(0,1){30}}
\put(115,45){$-6$}
\put(240,85){\line(0,1){30}}
\put(215,45){$-5$}
\put(340,85){\line(0,1){30}}
\put(315,45){$-4$}
\put(440,85){\line(0,1){30}}
\put(415,45){$-3$}
\put(540,85){\line(0,1){30}}
\put(515,45){$-2$}
\put(640,85){\line(0,1){30}}
\put(615,45){$-1$}
\put(840,85){\line(0,1){30}}
\put(830,45){$1$}
\put(940,85){\line(0,1){30}}
\put(930,45){$2$}
\put(1040,85){\line(0,1){30}}
\put(1030,45){$3$}
\put(1140,85){\line(0,1){30}}
\put(1130,45){$4$}
\put(1240,85){\line(0,1){30}}
\put(1230,45){$5$}
\put(1340,85){\line(0,1){30}}
\put(1330,45){$6$}
\end{picture}

\noindent
Since the function $t\mapsto F_\mu(t-i)U_tx$ is holomorphic in the stip
\hbox{$\{z\colon0<\Im{z}<2\}$}, we have
\[
2\mu\Int_{\Gamma}F_\mu(t-i)U_txdt=0
\]
and thus (\ref{tu}) takes the following form:
\[
(Q_\mu U_{2i}+2\mu Q_\mu U_i+\mu^2Q_\mu)x=\mu^2\Int_{\Gamma}F_\mu(t)U_txdt
=\mu^22\pi i\,{\rm
Res}_{\hspace{-5mm}\raisebox{-1.4mm}{${}_{t=i}$}}\,F_\mu(t)U_tx.
\]
The function $t\mapsto U_tx$ is holomorphic (and therefore continuous), so
that
\begin{eqnarray*}
{\rm
Res}_{\hspace{-5mm}\raisebox{-1.4mm}{${}_{t=i}$}}\,F_\mu(t)U_tx\hspace{.15cm}=&
\hspace{-1.04cm}\lim\limits_{t\rightarrow i}(t-i)F_\mu(t)U_tx&\\
=&\hspace{-.32cm}\lim\limits_{t\rightarrow
i}(t-i)F_\mu(t)\lim\limits_{t\rightarrow i}U_tx&\\
=&{\rm Res}_{\hspace{-5mm}\raisebox{-1.4mm}{${}_{t=i}$}}\,F_\mu(t) U_ix
=\frac{-i}{\mu^2 2\pi}U_ix.
\end{eqnarray*}
We have thus proved that
\begin{equation}
\label{podstwzor}
(Q_\mu U_{2i}+2\mu Q_\mu U_i+\mu^2Q_\mu)x=U_ix\where x\in X_\infty.
\end{equation}
Now if $x\in X$ is an arbitrary element one can define
\[
x_n=\sqrt{\frac{n}{\pi}}\Int_{-\infty}^{+\infty}U_txe^{-nt^2}dt
\]
and it is easily seen that for all $n\in\NN$ we have $x_n\in X_\infty$ and
$x_n\raisebox{1.2ex}{\hspace{.06cm}$\scriptscriptstyle{\sigma(X,\F)}$
\hspace{-.24cm}}
\raisebox{-0.7ex}{\hspace{-.52cm}$\scriptscriptstyle{n\rightarrow\infty}$\hspace
{.07cm}}
\mbox{\hspace{-.93cm}$-\!\!-\!\!\!\longrightarrow$}x$.
It is also easily verified with help of Proposition \ref{komut} that if
$x\in D(U_i)$ then
$U_ix_n\raisebox{1.2ex}{\hspace{.06cm}$\scriptscriptstyle{\sigma(X,\F)}$\hspace{
-.24cm}}
\raisebox{-0.7ex}{\hspace{-.52cm}$\scriptscriptstyle{n\rightarrow\infty}$\hspace
{.07cm}}
\mbox{\hspace{-.93cm}$-\!\!-\!\!\!\longrightarrow$}U_ix$.
This shows that $X_\infty$ is a core for $U_i$.

\noindent
On the Banach space $X\times X$ define a bounded operator
\[
\widetilde{R_\mu}=\begin{pmatrix}-Q_\mu+\frac{1}{\mu}I&-\frac{1}{\mu}Q_\mu\\
                                 \mu Q_\mu& Q_\mu\end{pmatrix}.
\]
Since $Q_\mu$ is $\sigma(X,\F)$-continuous so is $\widetilde{R_\mu}$. We shall
see that $\widetilde{R_\mu}$ leaves the \hbox{$\sigma(X,\F)$-closed} subspace
$\G(U_i)\subset X\times X$ invariant. Indeed: let $x\in X_\infty$. We have
\[
\widetilde{R_\mu}\begin{pmatrix}x\\ U_ix\end{pmatrix}
=\begin{pmatrix}-Q_\mu+\frac{1}{\mu}I&-\frac{1}{\mu}Q_\mu\\
                                 \mu Q_\mu&Q_\mu\end{pmatrix}
\begin{pmatrix}x\\ U_ix\end{pmatrix}
=\begin{pmatrix}\bigl(-Q_\mu+\frac{1}{\mu}I-\frac{1}{\mu}Q_\mu U_i\bigr)x\\
                 \bigl(\mu Q_\mu+Q_\mu U_i\bigr)x\end{pmatrix}.
\]
Now using (\ref{podstwzor}) and Proposition \ref{komut} we obtain
\[
U_i\left(-Q_\mu+\frac{1}{\mu}I-\frac{1}{\mu}Q_\mu U_i\right)x=
\left(-Q_\mu U_i+\frac{1}{\mu}U_i-\frac{1}{\mu}Q_\mu U_{2i}\right)x
=\bigl(\mu Q_\mu+Q_\mu U_i\bigr)x
\]
which means that $\widetilde{R_\mu}\bsm x\\ U_ix\esm\in\G(U_i)$.

\noindent
Remembering that $X_\infty$ is a core for $U_i$ and that $\widetilde{R_\mu}$
is \hbox{$\sigma(X,\F)$-continuous} one easily sees that
$\widetilde{R_\mu}\bigl(\G(U_i)\bigr)\subset\G(U_i)$. We can therefore consider
a new
\hbox{$\sigma(X,\F)$-continuous} operator
\[
R_\mu=\widetilde{R_\mu}|_{\G(U_i)}.
\]
Define
\[
\G_\infty=\left\{\begin{pmatrix} x\\y \end{pmatrix}\in
\G(U_i)\colon x,y\in X_\infty\right\}.
\]\sloppy
It is easily seen that $\G_\infty$ is a \hbox{$\sigma(X,\F)$-sequential} core
for $\Delta$ i. e. for any \hbox{$\bsm x\\ y\esm\in D(\Delta)$} there exists a
sequence
$\left\{\bsm x_n\\ y_n\esm\right\}_{n\in\NN}\subset\G_\infty$ such that
\begin{eqnarray}
&\begin{pmatrix}x_n\\ y_n\end{pmatrix}
\raisebox{1.2ex}{\hspace{.06cm}$\scriptscriptstyle{\sigma(X,\F)}$\hspace{-.24cm}
}
\raisebox{-0.7ex}{\hspace{-.52cm}$\scriptscriptstyle{n\rightarrow\infty}$\hspace
{.07cm}}
\mbox{\hspace{-.93cm}$-\!\!-\!\!\!\longrightarrow$}
\begin{pmatrix}x\\y\end{pmatrix},&\label{cwz1}\\
&\Delta\begin{pmatrix}x_n\\ y_n\end{pmatrix}
\raisebox{1.2ex}{\hspace{.06cm}$\scriptscriptstyle{\sigma(X,\F)}$\hspace{-.24cm}
}
\raisebox{-0.7ex}{\hspace{-.52cm}$\scriptscriptstyle{n\rightarrow\infty}$\hspace
{.07cm}}
\mbox{\hspace{-.93cm}$-\!\!-\!\!\!\longrightarrow$}
\Delta\begin{pmatrix}x\\y\end{pmatrix}.\label{cwz2}&
\end{eqnarray}

\begin{lem}
\label{coreinfty}
For any $\bsm x\\y\esm\in\G_\infty$ we have
\begin{eqnarray*}
&(\Delta+\mu I)R_\mu\bsm x\\y\esm=\bsm x\\y\esm,&\\
&R_\mu(\Delta+\mu I)\bsm x\\y\esm=\bsm x\\y\esm.
\end{eqnarray*}
\end{lem}
\dow
We compute
\begin{eqnarray*}
R_\mu\bigl(\Delta+\mu I\bigr)\begin{pmatrix} x\\U_ix\end{pmatrix}
\hspace{.22cm}=&\hspace{-3.36cm}\begin{pmatrix}
  -Q_\mu+\frac{1}{\mu}I&-\frac{1}{\mu}Q_\mu\\
  \mu Q_\mu&Q_\mu
  \end{pmatrix}
\begin{pmatrix}U_ix+\mu x\\U_{2i}x+\mu U_ix\end{pmatrix}&\\
=&\hspace{-1.19cm}\begin{pmatrix}
        \bigl(-Q_\mu U_i+\frac{1}{\mu}U_i-\mu Q_\mu+I-\frac{1}{\mu}Q_\mu
              U_{2i}-Q_\mu U_i\bigr)x\\
        \bigl(\mu Q_\mu U_i+\mu^2Q_\mu+Q_\mu U_{2i}+\mu Q_\mu U_i\bigr)x
  \end{pmatrix}&\\
=&\begin{pmatrix}
        \Bigl(\frac{1}{\mu}\bigl(U_i-Q_\mu U_{2i}-2\mu Q_\mu
U_i-\mu^2Q_\mu\bigr)+I\Bigr)x\\
        \bigl(Q_\mu U_{2i}+2\mu Q_\mu U_i+\mu^2Q_\mu\bigr)x
  \end{pmatrix}
\buildrel\ref{podstwzor}\over=\begin{pmatrix}x\\U_ix\end{pmatrix}.
\end{eqnarray*}
The formula
\[
\bigl(\Delta+\mu I\bigr)R_\mu\begin{pmatrix} x\\U_ix\end{pmatrix}
=\begin{pmatrix} x\\U_ix\end{pmatrix}
\]
may be derived analogously with a prior use of Proposition \ref{komut}.
\qed

\begin{lem}
\label{dziedzina}
We have $R_\mu\bigl(\G(U_i)\bigr)=D(\Delta)$.
\end{lem}
\dow
\underline{Step 1:} ``$R_\mu\bigl(\G(U_i)\bigr)\subset D(\Delta)$''. Take any
$\bsm x\\y\esm\in\G(U_i)$. then there is a sequence
$\{\bsm x_n\\y_n\esm\}_{n\in\NN}$ with
$\bsm x_n\\y_n\esm
\raisebox{1.2ex}{\hspace{.06cm}$\scriptscriptstyle{\sigma(X,\F)}$\hspace{-.24cm}
}
\raisebox{-0.7ex}{\hspace{-.52cm}$\scriptscriptstyle{n\rightarrow\infty}$\hspace
{.07cm}}
\mbox{\hspace{-.93cm}$-\!\!-\!\!\!\longrightarrow$}
\bsm x\\y\esm$.
Furthermore by the \hbox{$\sigma(X,\F)$-continuity} of $R_\mu$
we have
$R_\mu\bsm x_n\\y_n\esm
\raisebox{1.2ex}{\hspace{.06cm}$\scriptscriptstyle{\sigma(X,\F)}$\hspace{-.24cm}
}
\raisebox{-0.7ex}{\hspace{-.52cm}$\scriptscriptstyle{n\rightarrow\infty}$\hspace
{.07cm}}
\mbox{\hspace{-.93cm}$-\!\!-\!\!\!\longrightarrow$}
R_\mu\bsm x\\y\esm$ which combined with Lemma \ref{coreinfty} yields
\[
\bigl(\Delta+\mu I\bigr)R_\mu\begin{pmatrix}x_n\\y_n\end{pmatrix}
=\begin{pmatrix}x_n\\y_n\end{pmatrix}
\raisebox{1.2ex}{\hspace{.06cm}$\scriptscriptstyle{\sigma(X,\F)}$\hspace{-.24cm}
}
\raisebox{-0.7ex}{\hspace{-.52cm}$\scriptscriptstyle{n\rightarrow\infty}$\hspace
{.07cm}}
\mbox{\hspace{-.93cm}$-\!\!-\!\!\!\longrightarrow$}
\begin{pmatrix}x\\y\end{pmatrix}.
\]
Now since the operator $(\Delta+\mu I)$ is $\sigma(X,\F)$-closed we have
\begin{eqnarray}
&R_\mu\bsm x\\y\esm\in D(\Delta+\mu I)=D(\Delta),&\nonumber\\
&(\Delta+\mu I)R_\mu\bsm x\\y\esm=\bsm x\\y\esm.&\label{rezolw1}
\end{eqnarray}
In particular $R_\mu\bigl(\G(U_i)\bigr)\subset D(\Delta)$.

\noindent\sloppy
\underline{Step 2:} ``$R_\mu\bigl(\G(U_i)\bigr)\supset D(\Delta)$''. Let
$\bsm x\\ y\esm\in D(\Delta)$. Take a sequence
$\left\{\bsm x_n\\ y_n\esm\right\}_{n\in\NN}$ in the set $\G_\infty$ such that
formulae (\ref{cwz1}) and (\ref{cwz2}) hold. By Lemma \ref{coreinfty} we have
\begin{equation}
\label{RHSLHS}
R_\mu\bigl(\Delta+\mu I\bigr)\begin{pmatrix} x_n\\ y_n\end{pmatrix}
=\begin{pmatrix} x_n\\ y_n\end{pmatrix}\where n\in\NN,
\end{equation}
Taking the limit of both sides of (\ref{RHSLHS}) with $n\rightarrow\infty$
one obtains
\begin{equation}
\label{rezolweq}
R_\mu\bigl(\Delta+\mu I\bigr)\begin{pmatrix} x\\ y\end{pmatrix}
=\begin{pmatrix} x\\ y\end{pmatrix}.
\end{equation}
In particular $\bsm x\\ y\esm\in R_\mu\bigl(\G(U_i)\bigr)$ and consequently
$R_\mu\bigl(\G(U_i)\bigr)\supset D(\Delta)$.
\qed

\noindent
Combining Lemma \ref{dziedzina} and formulae (\ref{rezolw1}) and
(\ref{rezolweq}) we get the following proposition:

\begin{stw}
\label{rezolw}
The operator $R_\mu$ is the resolvent of the operator $\Delta$ at
\hbox{$-\mu\in\CC\setminus\RR_-$} i. e.
\[
\boxed{
\mbox{for any }\begin{pmatrix}x\\y\end{pmatrix}\in D(\Delta)\mbox{ we have }
R_\mu\bigl(\Delta+\mu I\bigr)
\begin{pmatrix}x\\y\end{pmatrix}=\begin{pmatrix}x\\y\end{pmatrix}
}
\]
and
\[
\boxed
{
\begin{array}{c}
\mbox{for any }\begin{pmatrix}x\\y\end{pmatrix}\in\G(U_i)
\mbox{ we have }R_\mu\begin{pmatrix}x\\y\end{pmatrix}\in D(\Delta)\\
\mbox{ and }\bigl(\Delta+\mu I\bigr)R_\mu\begin{pmatrix}x\\y\end{pmatrix}
=\begin{pmatrix}x\\y\end{pmatrix}.
\end{array}
}
\]
\end{stw}

\begin{wn}
\label{wynik}
Let $(X,\F)$ be a dual pair of Banach spaces having the Krein property and
such that the pair $(\F,X)$ also has the Krein property.
Let $U$ be a \hbox{$\sigma(X,\F)$-continuous}, one--parameter group of
isometries in $B_\F(X)$ and let $U_i$ be its analytic generator. Denote by
$\Delta$ the ampliation of $U_i$ onto $\G(U_i)$:
\[
\Delta=\begin{pmatrix}U_i&0\\0&U_i\end{pmatrix}.
\]
Then
\[
\Sp{\Delta}\subset\RR_+.
\]
Moreover the resolvent of $\Delta$ at $-\mu\in\CC\setminus\RR_-$ is given by
\[
R(\Delta,-\mu)=R_\mu=\begin{pmatrix}
                                 -Q_\mu+\frac{1}{\mu}I&-\frac{1}{\mu}Q_\mu\\
                                 \mu Q_\mu& Q_\mu
                     \end{pmatrix},
\]
where the operator $Q_\mu\in B_\F(X)$ is given by {\rm (\ref{Q})}.
\end{wn}

\noindent
For an operator $T$ on  $X$ let ${\rm Sp_{p}}T$ denote the point spectrum of
$T$. We have:

\begin{wn}
With notation as in Corollary \ref{wynik} we have
${\rm Sp_{p}}U_i\subset\RR_+$.
\end{wn}

\section{The Cior\v{a}nescu-Zsid\'{o} formula}

In the paper \cite{cizs} Ioana Cior\v{a}nescu and L\'{a}szl\'{o} Zsid\'{o}
proved that if $U$ is a \hbox{$\sigma(X,\F)$-continuous},
one--parameter group of isometries in $B_\F(X)$ then for any $t\in\RR$ we can
find $U_t$ using solely the analytic generator of $U$. In other words the
analytic generator determines the group. Their formula looks as follows
\[
U_tx=\lim_{\begin{smallmatrix}\alpha\rightarrow it\\O<{\rm
Re\,}\alpha<1\end{smallmatrix}}
\frac{sin\pi\alpha}{\pi}\Int_0^{+\infty}\lambda^{\alpha-1}
(\lambda+U_i)^{-1}U_ixd\lambda\where x\in D(U_i).
\]
Using the results of the previous section one can derive a similar formula --
namely for any $\bsm x\\y\esm\in D(\Delta)$ (i. e. for any $x\in D(U_{2i})$)
we have:
\[
U_tx=\lim_{\begin{smallmatrix}{\rm Im\,}z>0\\z\rightarrow t\end{smallmatrix}}
\frac{\sin-i\pi z}{\pi}\Int_0^{+\infty}\mu^{-iz-1}\Pr\nolimits_1(\Delta+\mu
I)^{-1}\Delta
\bsm x\\y\esm d\mu,
\]
where $\Pr\nolimits_1$ is the projection of $X\times X$ onto the first
coordinate. In our approach the computational details seem to be somewhat less
tedious. We shall not derive this formula here since it is mainly a repetition
of the work contained in \cite{cizs}, but one should remark that the family of
inequalities (cf \cite{cizs} p. 345)
\[
\|\Pr\nolimits_1
(\Delta+\mu I)^{-1}\Delta\bsm x\\y\esm\|
=\|\mu Q_\mu x+Q_\mu U_ix\|\leq\mu^{-r}C_r
\]
(where $r$ is a parameter in $]0,1[$ and $C_r$ a suitable constant) needed
to prove the convergence of the Bochner integral
\[
\Int_0^{+\infty}\mu^{\alpha-1}\Pr\nolimits_1(\Delta+\mu I)^{-1}\Delta\bsm
x\\y\esm d\mu
\]
(where $\alpha$ is any complex number such that $0<{\rm Re\,}\alpha<1$) can be
obtained by a simple computation:
\begin{eqnarray*}
\mu Q_\mu x+Q_\mu U_ix\hspace{.22cm}=&\hspace{-2.8cm}\Int_{-\infty}^{+\infty}\mu
F_\mu(t)U_txdt+
\Int_{-\infty}^{+\infty}F_\mu(t)U_{t+i}xdt&\\
=&\hspace{-1.82cm}
\Int_{-\infty}^{+\infty}\mu
F_\mu(t)U_{t+i}xdt+\Int_{i-\infty}^{i+\infty}F_\mu(\tau-i)U_\tau
xd\tau&\\
=&\hspace{-1.89cm}
\Int_{ir-\infty}^{ir+\infty}\mu
F_\mu(z)U_zxdz+\Int_{ir-\infty}^{ir+\infty}F_\mu(z-i)U_zxdz&\\
=&\hspace{-.22cm}\Int_{ir-\infty}^{ir+\infty}\bigl(\mu
F_\mu(z)+F_\mu(z-i)\bigr)U_zxdz
\buildrel(\ref{wl2})\over=\Int_{ir-\infty}^{ir+\infty}\frac{i\mu^{iz}U_zx}
{e^{\pi z}-e^{-\pi z}}
dz&\\
=&\hspace{-.44cm}
\Int_{-\infty}^{+\infty}\frac{i\mu^{i(t+ir)}U_{t+ir}x}{e^{\pi(t+ir)}-e^{-\pi(t+i
r)}}dt
=\mu^{-r}\Int_{-\infty}^{+\infty}\frac{i\mu^{it}U_tU_{ir}x}{e^{\pi(t+ir)}-
e^{-\pi(t+ir)}}dt.
\end{eqnarray*}

\section{Acknowledgements}

The author is greatly indebted to professor S. L. Woronowicz whose help in this
work could not be overestimated. This paper is a part of the author's master
thesis written under professor Woronowicz's supervision. The author would also
like to thank professor L\'{a}szl\'{o} Zsid\'{o} for fruitful discussions and
much help on subjects related to analytic generators.

\end{document}